\documentclass[preprint,12pt]{elsarticle}
\usepackage{amsmath}
\usepackage{amssymb}
\usepackage{color}
\usepackage{subfigure}
\usepackage{framed,multirow}
\usepackage{algorithm}
\usepackage{algorithmic}
\usepackage{graphicx}
\usepackage{subfigure}
\usepackage{epstopdf}
\usepackage{rotating}
\graphicspath{{Figure/}}
\begin{document}
	\begin{frontmatter}
	\title{MIM: A deep mixed residual method for solving high-order partial differential equations}
	\author[label1,label2]{Liyao Lyu}
	\ead{lyuliyao@msu.edu}
	\author[label3]{Zhen Zhang}
	\ead{zhangzhen@tbs-info.com}
	\author[label1,label4]{Minxin Chen\corref{cor1}}
	\ead{chenminxin@suda.edu.cn}
	\author[label1,label4]{Jingrun Chen\corref{cor1}}
	\ead{jingrunchen@suda.edu.cn}	
	\address[label1]{School of Mathematical Sciences, Soochow University, Suzhou, 215006, China}
	\address[label2]{CW Chu College, Soochow University, Suzhou, 215006, China}
	\address[label3]{Nanjing TBS Information Technology Co. Ltd, Nanjing, 210000, China}
	\address[label4]{Mathematical Center for Interdisciplinary Research, Soochow University, Suzhou, 215006, China}
	\cortext[cor1]{Corresponding authors.}
	
	\begin{abstract}
		In recent years, a significant amount of attention has been paid to solve partial differential equations (PDEs) by deep learning. For example, deep Galerkin method (DGM) uses the PDE residual in the least-squares sense as the loss function and a deep neural network (DNN) to approximate the PDE solution. In this work, we propose a deep mixed residual method (MIM) to solve PDEs with high-order derivatives. Notable examples include Poisson equation, Monge-Amp\'{e}re equation, biharmonic equation, and Korteweg-de Vries equation. In MIM, we first rewrite a high-order PDE into a first-order system, very much in the same spirit as local discontinuous Galerkin method and mixed finite element method in classical numerical methods for PDEs. We then use the residual of first-order system in the least-squares sense as the loss function, which is in close connection with least-squares finite element method. For aforementioned classical numerical methods, the choice of trail and test functions is important for stability and accuracy issues in many cases. MIM shares this property when DNNs are employed to approximate unknowns functions in the first-order system. In one case, we use nearly the same DNN to approximate all unknown functions and in the other case, we use totally different DNNs for different unknown functions. Numerous results of MIM with different loss functions  and different choice of DNNs are given for four types of PDEs. In most cases, MIM provides better approximations (not only for high-derivatives of the PDE solution but also for the PDE solution itself) than DGM with nearly the same DNN and the same execution time, sometimes by more than one order of magnitude. When different DNNs are used, in many cases, MIM provides even better approximations than MIM with only one DNN, sometimes by more than one order of magnitude. Numerical observations also imply a successive improvement of approximation accuracy when the problem dimension increases and interesting connections between MIM and classical numerical methods. Therefore, we expect MIM to open up a possibly systematic way to understand and improve deep learning for solving PDEs from the perspective of classical numerical analysis.
	\end{abstract}
\end{frontmatter}
\section{Introduction}
\label{sec:Intro}
	Solving partial differential equations (PDEs) has been the most ubiquitous tool to simulate complicated phenomena in applied sciences and engineering problems. Classical numerical methods include finite difference method \cite{leveque_finite_2007}, finite element method (FEM) \cite{Elman2014},
	discontinuous Galerkin method \cite{Cockburn2000}, and spectral method \cite{shen2011spectral}, which are typically designed for low dimensional PDEs and are well understood in terms of stability and accuracy. However, there are high dimensional PDEs such as Schr\"odinger equation in the quantum many-body problem \cite{dirac1981principles}, Hamilton-Jacobi-Bellman equation in stochastic optimal control \cite{bardi2008optimal}, and nonlinear Black-Scholes equation for pricing financial derivatives \cite{hull2009options}. Solving these equations is far out of the capability of classical numerical methods due to the curse of dimensionality, i.e., the number of unknowns grows exponentially fast as the dimension increases.
	
	Until very recently, deep-learning based methods have been developed to solving these high-dimensional PDEs; see \cite{Weinan2017349,Carleo2017,E2018Mar,Jiequn2018,raissi2018deep,sirignano2018dgm, Hutzenthaler2019Jan,raissi2019physics-informed,Beck20191563,GonzlezCervera2019,fan2019multiscale,khoo2019solving, Beck2020Mar,wang2020deep,zang2020weak,discacciati2020controlling} for examples. Typically, there are three main ingredients (stages) of a deep-learning method for solving PDEs: (1) modeling: the loss (objective) function to be optimized; (2) architecture: the deep neural network (DNN) for function approximation; (3) optimization: the optimal set of parameters in the DNN which minimizes the loss function. By design, the number of parameters in DNNs grows at most polynomially in terms of dimension. Meanwhile, possibly high-dimensional integrals in the loss function are approximated by Monte-Carlo method. Therefore, by design, deep learning overcomes the curse of dimensionality. In practice, deep learning performs well for Schr\"odinger equation \cite{Carleo2017,Han2019}, Hamilton-Jacobi-Bellman equation \cite{Jiequn2018,Weinan2017349}, and nonlinear Black-Scholes equation \cite{Beck20191563,GonzlezCervera2019}.
	
	Typically, deep learning solves a PDE in the following way. For the given PDE, the loss function is modeled as the equation residual in the least-squares sense \cite{sirignano2018dgm} or the variational form if exists \cite{E2018Mar}. ResNet is often used as the network architecture \cite{He2015}, which was tested to overcome the notorious problem of vanishing/exploding gradient. Afterwards, stochastic gradient descent method is used to find the optimal set of parameters in ResNet which minimizes the loss function. ResNet with the optimal set of parameters gives an approximation of the PDE solution.
	
	In this work, we propose a deep mixed residual method (MIM) for solving high-order PDEs. In the modeling stage, by rewriting a given PDE into a first-order system, we obtain a larger problem in the sense that both the PDE solution and its high-order derivatives are unknown functions to be approximated. This has analogs in classical numerical methods, such as local discontinuous Galerkin method \cite{Cockburn2000} and mixed finite element method \cite{Boffi2013}.
	Compared to DGM, there are two more degrees of freedom in MIM:
	\begin{itemize}
		\item In the loss function stage, one can choose different high-order derivatives into the set of unknown functions. Take biharmonic equation as an example. The set of unknown functions can include the PDE solution and its derivatives up to the third order, or only contain the PDE solution and its second-order derivatives, and both choices have analogs in discontinuous Galerkin method \cite{Yan2002Dec, cockburn2009a}. We then write the loss function as the sum of equation residuals in the least-squares sense, very much in the same spirit as the least-squares finite element method \cite{Bochev2015}.
		\item In the architecture stage, one can choose the number of networks to approximate the set of unknown functions. In one case, one DNN is used to approximate the PDE solution and other DNNs are used to approximate its high-order derivatives; in the other case, the PDE solution and its derivatives share nearly the same DNN.
	\end{itemize}
	These two degrees of freedom allow MIM to produce better approximations over DGM in all examples, including Poisson equation, Monge-Amp\'{e}re equation, biharmonic equation, and Korteweg-de Vries (KdV) equation. In particular, MIM provides better approximations not only for the high-order derivatives but also for the PDE solution itself. It is worth mentioning that the usage of mixed residual in deep learning was first introduced for surrogate modeling and uncertainty quantification of a second-order elliptic equation \cite{Zhu2019Oct} and was later adopted in a deep domain decomposition method \cite{Li2019Dec}.
			
	The paper is organized as follows. In Section \ref{sec:mixed residual}, we introduce MIM and DGM (for comparison purpose). In Section \ref{sec:Numerical result}, numerical results for four types of high-order PDEs are provided. Conclusions and discussions are drawn in Section \ref{sec:conclusion}.
	
	\section{Deep mixed residual method}
	\label{sec:mixed residual}
	
	In this section, we introduce MIM and discuss its difference with DGM in terms of loss function and neural network structure.
	
	\subsection{Loss function}
	\label{subsection:Loss function}

	Consider a potentially time-dependent nonlinear PDE over a bounded domain $\Omega\subset\mathbb{R}^d$
   \begin{equation}
   \left\{
   	\begin{aligned}
   		\label{equ:gen time_dependent pde example}
   		&\partial_t u + \mathcal{L} u = 0 & (t,x)\in (0,T]\times\Omega,\\
  		 &u(0,x)= u_0(x) & x\in \Omega, \\
   		&u(t,x) = g(x) & (t,x)\in [0,T]\times \partial \Omega,
   	\end{aligned}
   	\right.
   \end{equation}
   where $\partial \Omega$ denotes the boundary of $\Omega$. In DGM, the loss function is defined as the PDE residual in the least-squares sense
   \begin{equation}
   \label{eqn:lossdgm}
	\begin{aligned}
		L(u) = \|\partial_t u + \mathcal{L} u\|^2_{2,[0,T]\times\Omega}
		+ \lambda_1\|u(0,x)-u_0\|^2_{2,\Omega}
		+ \lambda_2\|u - g\|^2_{2, [0,T]\times \partial \Omega},
	\end{aligned}
   \end{equation}
   where $\lambda_1$ and $\lambda_2$ are penalty parameters given \textit{a priori}. These three terms in \eqref{eqn:lossdgm} measure how well the approximate solution satisfies the PDE, the initial condition and the boundary condition, respectively.

   In the absence of temporal derivatives, \eqref{equ:gen time_dependent pde example} reduces to
   \begin{equation*}
   \left\{
			   	\begin{aligned}
			   		&\mathcal{L} u = 0 & x\in \Omega,\\
			   		& u(x)= g(x) & x\in \partial \Omega,
			   	\end{aligned}
			   	\right.
   \end{equation*}
   and the corresponding loss function in DGM becomes
   \begin{equation}\label{equ:loss function}
			   	\begin{aligned}
			   	L(u) = \|\mathcal{L}u\|^2_{2,\Omega} + \lambda \|u-g\|^2_{2,\partial \Omega}.
			   	\end{aligned}
   \end{equation}

   Table \ref{tbl:loss DGM} lists four PDEs with their corresponding loss functions in DGM and Table \ref{tbl:loss BC} lists different boundary conditions, the initial condition and their contributions to loss functions in DGM and MIM. More boundary conditions can be treated in this way. Interested readers may refer to \cite{chen2020bc} for details.
   \begin{table}\centering
   	\begin{tabular}{|l|l|l|}
   		\hline
   	 Equation & Explicit form &  Loss function $L(u)$ \\
   	 \hline
   	 Poisson &
   	 $
   	 -\Delta u = f(x)
   	 $ & $\|\Delta u + f(x) \|_{2,\Omega}^2$
   	 \\
   	 \hline
   	 Monge-Amp\'ere &
   	 $\det(\nabla^2 u) = f(x)$ & $\|\det(\nabla^2u) - f(x)\|_{2,\Omega}^2 $ \\
   	 \hline
   	 Biharmonic &
   	 $
   	 -\Delta^2 u = f(u,x)
   	 $&
   	 $\|
   	 \Delta^2 u + f(u,x)
   	 \|_{2,\Omega}^2
   	 $\\
   	 \hline
   	 KdV &
   	 $
   	 u_t + \sum_{i=1}^{d}u_{x_ix_ix_i} = f(x) $
   	 &
   	 $
   	 \|
   	 u_t + \sum_{i=1}^{d}u_{x_ix_ix_i} - f(x)
   	 \|_{2,\Omega}^2
   	 $
   	 \\
   	 \hline
   	 \end{tabular}
    \caption{Loss functions for four types of PDEs in the deep Galerkin method.}
	\label{tbl:loss DGM}
	\end{table}
	\begin{table}
		\centering\begin{tabular}{|l|l|l|}\hline
			Condition & Explicit form & Contribution to the loss function\\
			\hline
			Dirichlet &
			$u(x) = g$ & $\|u-g\|_{2,[0,T]\times\partial\Omega}^2$\\ \hline
			Neumann &
			$\frac{\partial u}{\partial n} = g$ & $\|\frac{\partial u}{\partial n} - g\|_{2,[0,T]\times\partial\Omega}^2$ or $\|p - g\|_{2,[0,T]\times\partial\Omega}^2$\\ \hline
			Initial &  $u(0,x) = u_0(x)$ & $\|u-u_0\|^2_{2,\Omega}$\\
			\hline
		\end{tabular}
	\caption{Contributions to the loss function for the initial condition and different types of boundary conditions used in the deep Galerkin method and the deep mixed residual method.}
	\label{tbl:loss BC}
	\end{table}
	
	In MIM, we first rewrite high-order derivatives into low-order ones using auxiliary variables. For notational convenience, auxiliary variables $p,q,w$  represent
	\begin{equation}\label{eqn:auxiliary}
		\begin{aligned}
		p &= \nabla u, \\
		q &= \nabla\cdot p = \Delta u, \\
		w &= \nabla q = \nabla (\Delta u).
		\end{aligned}
	\end{equation}
	For KdV equation, we have $q=\mathrm{diag}(\nabla p)$ instead of the second formula in \eqref{eqn:auxiliary}.
	With these auxiliary variables, we define loss functions for four types of PDEs in Table \ref{tbl:loss MIM}. Since one can choose a subset of high-order derivatives into the set of unknown functions, there are more than one loss function in MIM. For biharmonic equation, there are two commonly used sets of auxiliary variables in local discontinuous Galerkin method and weak Galerkin finite element method: one with all high-order derivatives \cite{Yan2002Dec} and the other with part of high-order derivatives \cite{cockburn2009a,Mu2015Sep}. Correspondingly, if all high-order derivatives are used, we denote MIM by MIM$_{a}$, and if only part of high-order derivatives are used, we denote MIM by MIM$_{p}$. In Section \ref{subsec:neural network}, we will discuss how to equip different loss functions with different DNNs. In short, if only one DNN is used to approximate the PDE solution and its derivatives, we denote MIM by MIM$^1$, and if multiple DNNs are used, we denote MIM by MIM$^2$. In Section \ref{sec:Numerical result}, different loss functions listed in Table \ref{tbl:loss DGM}, Table \ref{tbl:loss BC} and Table \ref{tbl:loss MIM} will be tested and discussed. By default, all the penalty parameters are set to be $1$.
	\begin{table}[htbp]\centering
		\resizebox{\textwidth}{!}{
		\begin{tabular}{|l|l|l|}
			\hline
			Equation & Explicit form &  Loss function $L(u,p,q,w)$ \\
			\hline
			Poisson &
			$
			-\Delta u = f(u,x)
			$
			& $\|p - \nabla u \|_{2,\Omega}^2 + \|\nabla \cdot p + f(u,x) \|_{2,\Omega}^2$
			\\
			\hline
			Monge-Amp\'ere &
			$\det(\nabla ^2 u) = f$
			& $\|p - \nabla u \|_{2,\Omega}^2+\|\det(\nabla p) - f\|_{2,\Omega}^2 $\\
			\hline
			\multirow{3}*{Biharmonic} &
			\multirow{3}*{
			$
			-\Delta^2 u = f(u,x)
			$}
			&$\|p - \nabla u \|_{2,\Omega}^2 +\|q - \nabla \cdot p\|_{2,\Omega}^2 $
			\\
			~ & ~ & $+ \|w - \nabla q\|_{2,\Omega}^2 + \|\nabla \cdot w + f\|_{2,\Omega}^2$\\
			\cline{3-3}
			~ & ~ & $\|q - \Delta u \|_{2,\Omega}^2 + \|\Delta q + f\|_{2,\Omega}^2$
			\\
			\hline
			\multirow{2}*{KdV} & \multirow{2}*{$u_t + \sum_{i=1}^{d}u_{x_ix_ix_i} = f(x)$}
			&
			$\|p - \nabla u\|_{2,[0,T]\times\Omega}^2 + \|q - \mathrm{diag}(\nabla p) \|_{2,[0,T]\times\Omega}^2 $\\
			~ & ~ & $+  \|u_t + \nabla\cdot q - f(x)\|_{2,[0,T]\times\Omega}^2 $ \\
			\hline
		\end{tabular}}
		\caption{Loss functions in the deep mixed residual method for four types of equations. Two different loss functions for biharmonic equation are denoted by MIM$_a$ and MIM$_p$, in which all high-order derivatives or part of high-order derivatives are included, respectively.}
		\label{tbl:loss MIM}
	\end{table}

	\subsection{Neural network architecture}
	\label{subsec:neural network}
	
\begin{figure}[ht]
		\centering
		\includegraphics[width=0.3\textwidth]{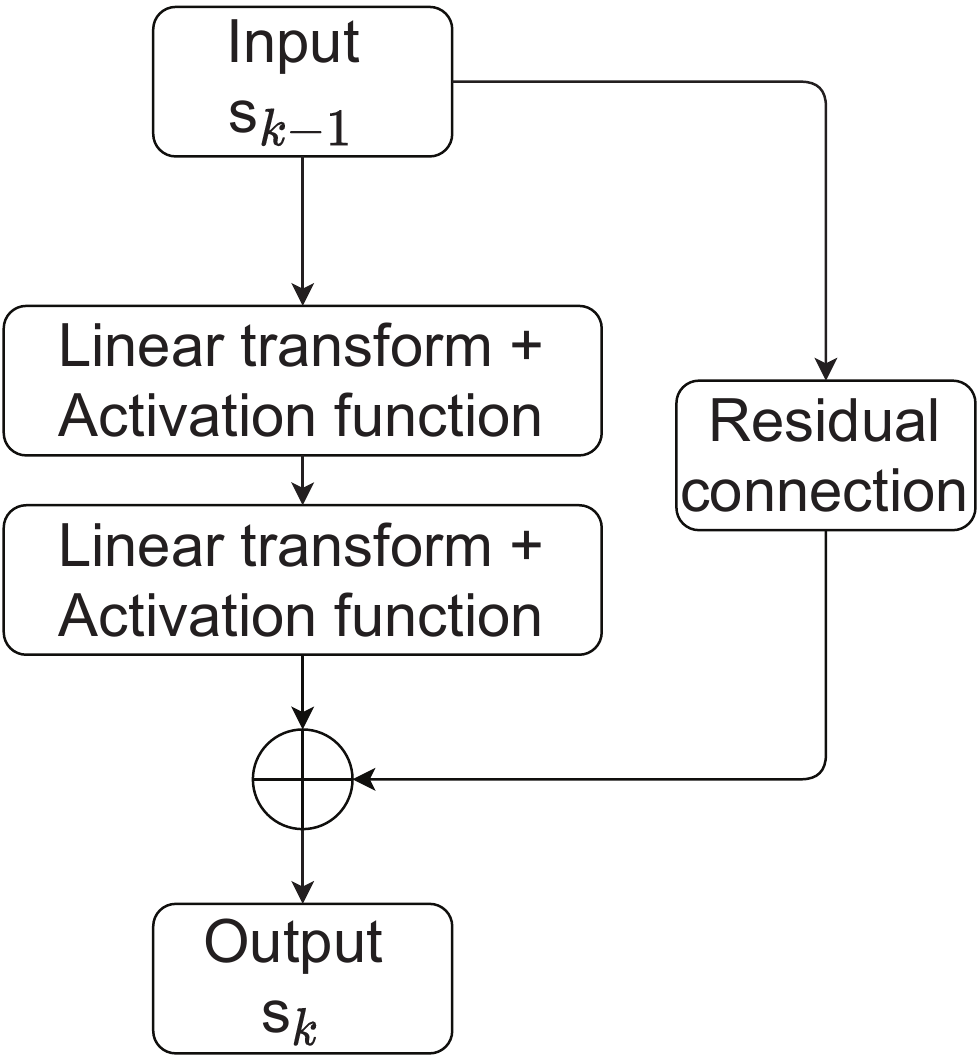}
		\caption{One block of ResNet. A deep neural network contains a sequence of blocks, each of which consists of two fully-connected layers and one shortcut connection.}
		\label{fig:resnet}
	\end{figure}

	ResNet \cite{He2015} is used to approximate the PDE solution and its high-order derivatives. It consists of $m$ blocks in the following form
	\begin{equation}\label{equ:resnet}
		s_k =\sigma(W_{2,k}\sigma(W_{1,k}s_{k-1}+b_{1,k})+b_{2,k}) +s_{k-1}, \quad k=1,2,\cdots,m.
	\end{equation}
	Here $s_k, b_{1,k}, b_{2,k} \in \mathbb{R}^n$, $W_{1,k}, W_{2,k}\in \mathbb{R}^{n\times n}$. $m$ is the depth of network, $n$ is the width of network, and $\sigma$ is the (scalar) activation function. Explicit formulas of activation functions used in this work are given in Table \ref{tbl:activaction function}. The last term on the right-hand side of \eqref{equ:resnet} is called the shortcut connection or residual connection. Each block has two linear transforms, two activation functions, and one shortcut; see Figure \ref{fig:resnet} for demonstration. Such a structure can automatically solve the notorious problem of vanishing/exploding gradient \cite{DBLP:journals/corr/HeZRS15}.
	\begin{table}
	  \centering
	  \begin{tabular}{|l|l|}
	    \hline
	    Activation function & Formula \\
	    \hline
	    Square & $x^2$ \\
	    \hline
	    ReLU & $\max\{x,0\}$\\
	    \hline
		ReQU & $(\max\{x,0\})^2$\\
		\hline
		ReCU & $(\max\{x,0\})^3$\\
	    \hline
	  \end{tabular}
	  \caption{Activation functions used in numerical tests.}\label{tbl:activaction function}
	\end{table}

\begin{figure}
		\centering
		\includegraphics[width=\textwidth]{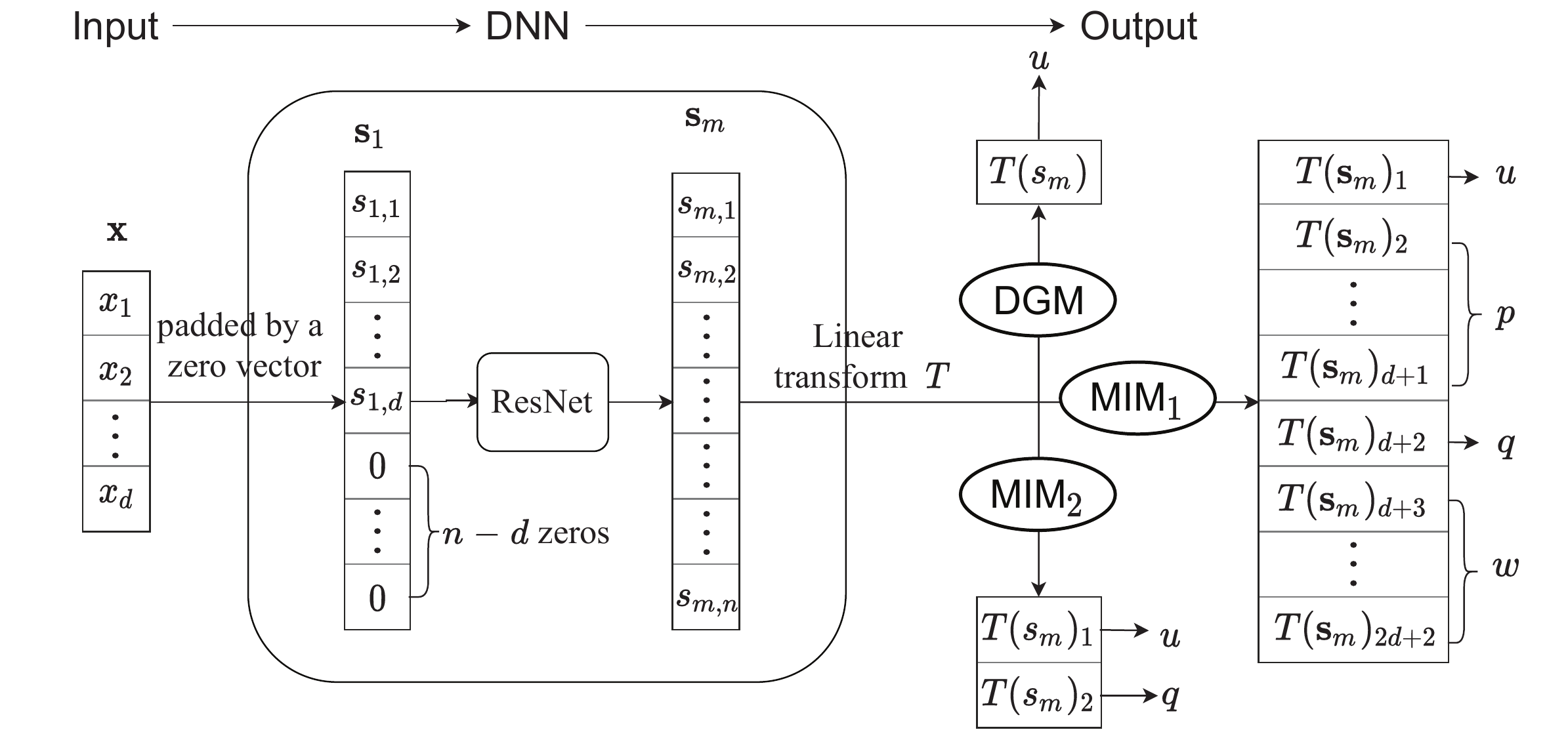}
		\caption{Network structures for biharmonic equation with deep Galerkin method and deep mixed residual method. DGM only approximates solution $u$. MIM$^1_p$ approximate solution $u$ and $\Delta u$. MIM$^1_a$ approximates solution $u$ and all of its derivatives used in the equation $\nabla u, \Delta u , \nabla (\Delta u)$. MIM$^2_a$ uses four networks to approximate $ u , \nabla u, \Delta u , \nabla (\Delta u)$ and MIM$^2_p$ uses two networks to approximate $u,\Delta u$. Each network has a similar structure with different output dimensions.}
		\label{fig:network for biharmonic}
	\end{figure}

	Since $x$ is in $\mathbb{R}^d$ rather than $\mathbb{R}^n$, we can pad $x$ by a zero vector to get the network input $s_0$. A linear transform can be used as well without much difference. Meanwhile, $s_m$ has $n$ outputs which cannot be directly used for the PDE solution and its derivatives employed in the loss function. Therefore, a linear transform $T$ is applied to $s_m$ to transform it into a suitable dimension. Let $\{\theta\}$ be the whole set of parameters  which include parameters in ResNet ($\left\{W_{1,k},b_{1,k},W_{2,k},b_{2,k}\right\}_{k=1}^m$) and parameters in the linear transform $T$. 
	Note that the output dimension in MIM depends on both the PDE problem and the mixed residual loss. We illustrate network structures for biharmonic equation as an example in Figure \ref{fig:network for biharmonic}.
	From Figure \ref{fig:network for biharmonic}, we see that DGM has only $1$ output, MIM$^1_a$ has $2d+2$ outputs, and MIM$^1_{p}$ has $2$ outputs.
	In Figure \ref{fig:neural network for poisson equation MIM}, we illustrate networks structures of MIM$^1$ and MIM$^2$ for Poisson equation.
	In MIM$^2$, two DNNs are used: one to approximate the solution and the other one to approximate its derivatives. It is clear from Figure \ref{fig:network for biharmonic} that network structures in DGM and MIM$^1$ only differ in the output layer and thus they have comparable numbers of parameters to be optimized. To be precise, we  calculate their numbers of parameters in Table \ref{tbl:param num}, from which one can see the number of parameters in DGM and MIM$^1$ is close. The number of parameters in MIM$^2$ is nearly double for Poisson equation, Monge-Amp\'ere equation and biharmonic equation (MIM$^2_p$), tripled for KdV equation, and quadrupled for biharmonic equation (MIM$^2_a$), respectively. In Section \ref{sec:Numerical result}, from numerical results, we observe a better performance of MIM$^1$ for all four equations, not only for derivatives of the PDE solution, but also for the solution itself.
	\begin{table}
		\centering\begin{tabular}{|l|l|l|}
			\hline
			 Method & Equation & Size of the parameter set\\
			 \hline
			DGM & Four equations &  $(2m-1)n^2 + (2m+d+1)n +1$\\
			\hline
			\multirow{5}*{MIM$^1$} & Poisson &\multirow{2}*{$(2m-1)n^2 + (2m+2d+1)n +d + 1$} \\
			\cline{2-2}
			~ & Monge-Amp\'ere & ~ \\
			\cline{2-3}
			~ & Biharmonic (MIM$^1_a$) & $(2m-1)n^2 + (2m+3d+2)n +2d + 2$  \\
			\cline{2-3}
			~ & Biharmonic (MIM$^1_p$) & $(2m-1)n^2 + (2m+d+2)n + 2$ \\
			\cline{2-3}
			~ & KdV & $(2m-1)n^2 + (2m+3d+1)n +2d + 1$ \\
			\hline
			\multirow{5}*{MIM$^2$} & Poisson &\multirow{2}*{$(4m-2)n^2 + (4m+3d+1)n +d + 1$} \\
			\cline{2-2}
			~ & Monge-Amp\'ere & ~ \\
			\cline{2-3}
			~ & Biharmonic (MIM$^2_a$) & $(8m-4)n^2 + (8m+6d+2)n +2d + 2$  \\
			\cline{2-3}
			~ & Biharmonic (MIM$^2_p$) & $(4m-2)n^2 + (4m+2d+2)n + 2$ \\
			\cline{2-3}
			~ & KdV & $(6m-3)n^2 + (6m+5d+1)n +2d + 1$ \\
			\hline
		\end{tabular}
	\caption{Number of parameters for different network structures used for different equations and different loss functions. $n$, $m$, and $d$ are the network width, the network depth, and the problem dimension, respectively. It is observed that the number of parameters in DGM and MIM$^1$ is close, and the number of parameters in MIM$^2$ is nearly double for Poisson equation, Monge-Amp\'ere equation and biharmonic equation (MIM$^2_p$), tripled for KdV equation, and quadrupled for biharmonic equation (MIM$^2_a$), respectively.}
	\label{tbl:param num}
	\end{table}
	
	\subsection{Stochastic Gradient Descent}
	
	For completeness, we also briefly introduce stochastic gradient descent method. For the loss function defined in \eqref{equ:loss function}, we generate two sets of points uniformly distributed over $\Omega$ and $\partial\Omega$: $\{\mathbf{x}_i\}^{N}_{i=1}$ in $\Omega$ and $\{\mathbf{\hat{x}}_j\}^{M}_{j=1}$ on $\partial \Omega$.
	\begin{equation}
		\theta^{k+1} = \theta^{k} - \alpha \nabla_{\theta} \frac{|\Omega|}{N} \sum_{i=1}^N [\mathcal{L}u_{\theta}(\mathbf{x}_i;\theta^{k})]^2 + \lambda \alpha \nabla_{\theta} \frac{|\partial \Omega|}{M} \sum_{j=1}^{M}  [u_{\theta}(\mathbf{\hat{x}}_j;\theta^{k})-g(\mathbf{\hat{x}}_j)]^2,
	\end{equation}
	where $\alpha$ is the learning rate chosen to be $1e-3$ here. $|\Omega|$ and $|\partial \Omega|$ are measures of $\Omega$ and $\partial\Omega$, respectively. $u_{\theta}$ is the DNN approximation of PDE solution parameterized by $\{\theta\}$. Sampling points $\{\mathbf{x}_i\}^{N}_{i=1}$ and $\{\mathbf{\hat{x}}_j\}^{M}_{j=1}$ are updated at each iteration. In implementation, we use ADAM optimizer \cite{kingma2015adam} and automatic differentiation \cite{Paszke2017Oct} for derivatives in PyTorch.
	
	\section{Numerical Result}
	\label{sec:Numerical result}

	In this section, we show numerical results of MIM for four types of equations. We use relative $L^2$ errors of $u$, $\nabla u$, $\Delta u$, and $\nabla (\Delta u)$ defined in Table \ref{tbl:R2error} for comparison. In all figures, relative $L^2$ errors are in $\log _{10}$ scale.
	\begin{table}
		\centering
		\begin{tabular}{|c|c|c|}\hline
			Quantity  & DGM  & MIM\\
			 \hline
			 $u$ & $\frac{\int_{\Omega} ( u_{\theta} - u ) ^2 dx}{\int_{\Omega} u ^2 dx}$ & $\frac{\int_{\Omega} (u_{\theta} - u ) ^2 dx}{\int_{\Omega} u ^2 dx}$ \\
			 \hline
			 $\nabla u $ & $\frac{\int_{\Omega} (\nabla u_{\theta} - \nabla u ) ^2 dx}{\int_{\Omega} (\nabla u) ^2 dx}$ &
			 $\frac{\int_{\Omega} (p_{\theta} - \nabla u ) ^2 dx}{\int_{\Omega} (\nabla u) ^2 dx}$
			  \\
			 \hline
			 $\Delta u $ & $\frac{\int_{\Omega} (\Delta u_{\theta} - \Delta u ) ^2 dx}{\int_{\Omega} (\Delta u) ^2 dx}$ &
			 $\frac{\int_{\Omega} (q_{\theta} - \Delta u ) ^2 dx}{\int_{\Omega} (\Delta u) ^2 dx}$\\
			 \hline
			 $\nabla \Delta u $ & $\frac{\int_{\Omega} \left(\nabla (\Delta u_{\theta}) - \nabla (\Delta u)\right) ^2 dx}{\int_{\Omega} \left(\nabla (\Delta u)\right) ^2 dx}$ &
			 $\frac{\int_{\Omega} \left(w_{\theta} - \nabla (\Delta u) \right) ^2 dx}{\int_{\Omega} \left(\nabla (\Delta u)\right) ^2 dx}$\\\hline
			 $\mathrm{diag}(\nabla^2 u) $ & $\frac{\int_{\Omega} (\mathrm{diag}(\nabla^2 u_{\theta} )- \mathrm{diag}(\nabla^2 u) ) ^2 dx}{\int_{\Omega} \left(\mathrm{diag}(\nabla^2 u)\right) ^2 dx}$ &
			 $\frac{\int_{\Omega} (q_{\theta} - \mathrm{diag}(\nabla^2 u) ) ^2 dx}{\int_{\Omega} \left(\mathrm{diag}(\nabla^2 u)\right) ^2 dx}$\\
			 \hline
		\end{tabular}
				 \caption{Relative $L^2$ errors used in deep Galerkin method and deep mixed residual method.}
			\label{tbl:R2error}
	\end{table}
	\subsection{Poisson Equation}
	Consider the following Neumann problem
	\begin{equation}\label{eq:case_Neumann Boundary condition}
	\left\{
	\begin{aligned}
	&-\Delta u + \pi^2 u = 2 \pi^2 \sum_{k=1}^{d} \cos(\pi x_k) & x\in \Omega=[0,1]^d\\
	&\frac{\partial u}{\partial n} = 0 & x\in \partial \Omega
	\end{aligned}\right.
	\end{equation}
	with the exact solution $u(x)=\sum_{k=1}^{d}\cos(\pi x_k)$. The neural network structure in DGM is the same as that for biharmonic equation shown in Figure \ref{fig:network for biharmonic}. Following Table \ref{tbl:loss DGM} and Table \ref{tbl:loss BC}, we use the loss function for \eqref{eq:case_Neumann Boundary condition}
	\begin{equation}
	\begin{aligned}
	L(u) = &\|-\Delta u + \pi^2 u - 2 \pi^2 \sum_{k=1}^{d} \cos(\pi x_k)\|_{2,\Omega}^2 +\lambda \|\frac{\partial u}{\partial n} \|_{2,\partial \Omega}^2.
	\end{aligned}
	\end{equation}
	
\begin{figure}
		\centering
		\subfigure[MIM$^1$: one network to approximate the PDE solution and its derivatives.]{
		\includegraphics[width=\textwidth]{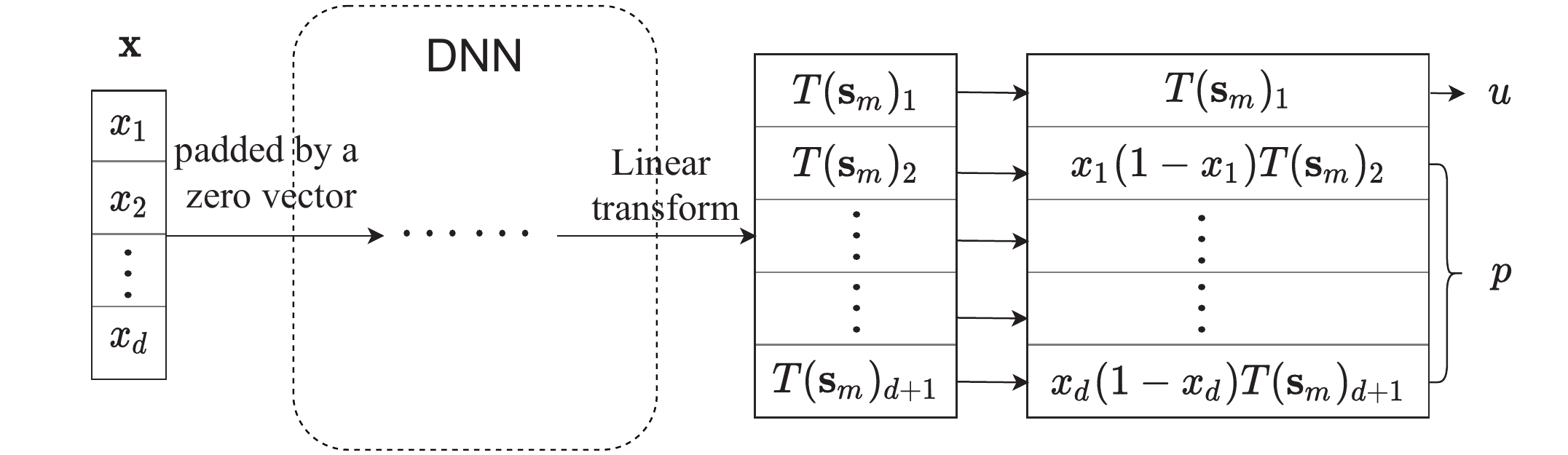}}
		\subfigure[MIM$^2$: multiple networks to appriximate the PDE solution and its derivatives.]{
		\includegraphics[width=\textwidth]{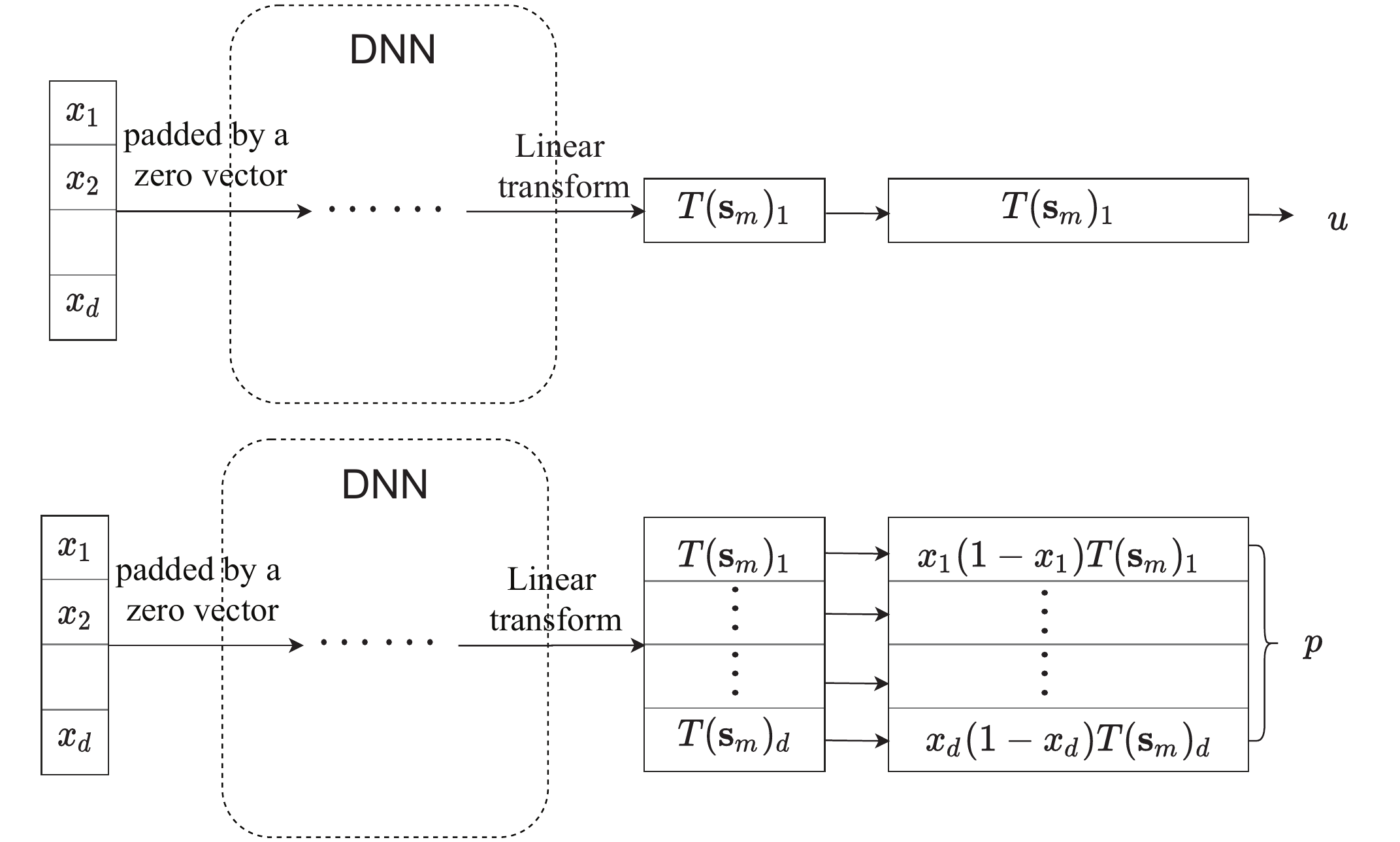}}
		\caption{Detailed network structures of MIM$^1$ and MIM$^2$ to solve Poisson equation. DNN part is the same as that in Figure \ref{fig:network for biharmonic}. $x_i(1-x_i)$ are multipliers which make MIM$^1$ and MIM$^2$ satisfy the exact Neumann boundary condition.}
		\label{fig:neural network for poisson equation MIM}
	\end{figure}

	Since both $u$ and $p$ are explicitly used, one more advantage of MIM is the enforcement of boundary conditions. For \eqref{eq:case_Neumann Boundary condition}, we multiply $p_i,\;i=1,\cdots,d$ by $x_i(1-x_i)$ to satisfy the Neumann boundary condition automatically; see Figure \ref{fig:neural network for poisson equation MIM}. DGM only has $u$ as its unknown function, and thus it is unclear that how the exact Neumann boundary condition can be imposed.
	Therefore, for DNNs in Figure \ref{fig:neural network for poisson equation MIM}, the loss function in MIM can be simplified as
	\begin{equation}
	L(u,p) = \|p- \nabla u\|^2_{2,\Omega} + \|-\nabla\cdot p+\pi^2u-2\pi^2\sum_{k=1}^{d}\cos(\pi x_k)\|^2_{2,\Omega}.
	\end{equation}
	We emphasize that Dirichlet boundary condition can be exactly imposed in DGM \cite{berg_unified_2018} and no penalty term is needed. For Neumann boundary condition, mixed boundary condition, and Robin boundary condition, however, it is difficult to build up a DNN representation which satisfies the exact boundary condition. Building up a DNN approximation which satisfies the exact boundary condition can have a couple of advantages \cite{chen2020bc}: 1) make ease of the training process by avoiding unnecessary divergence; 2) improve the approximation accuracy; 3) save the execution time. In MIM, however, we have the direct access to both $u$ and $p$. Therefore, all these boundary conditions can be imposed exactly in principle. This will be presented in a subsequent work \cite{lyu2020}.
	
\begin{figure}[ht]
		\centering
%
		\subfigure[4D]{
			\includegraphics[width=0.3\textwidth]{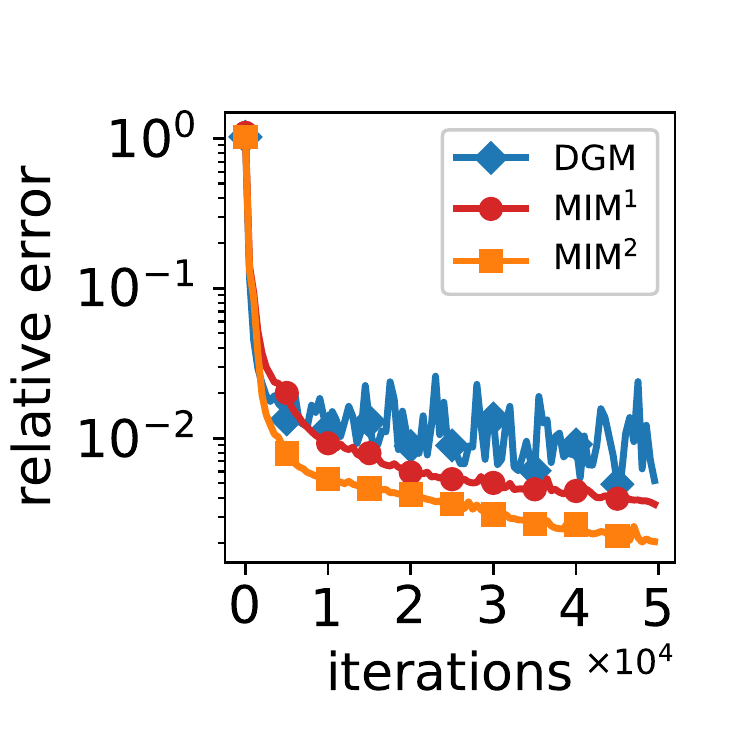}
		}
		\subfigure[8D]{
			\includegraphics[width=0.3\textwidth]{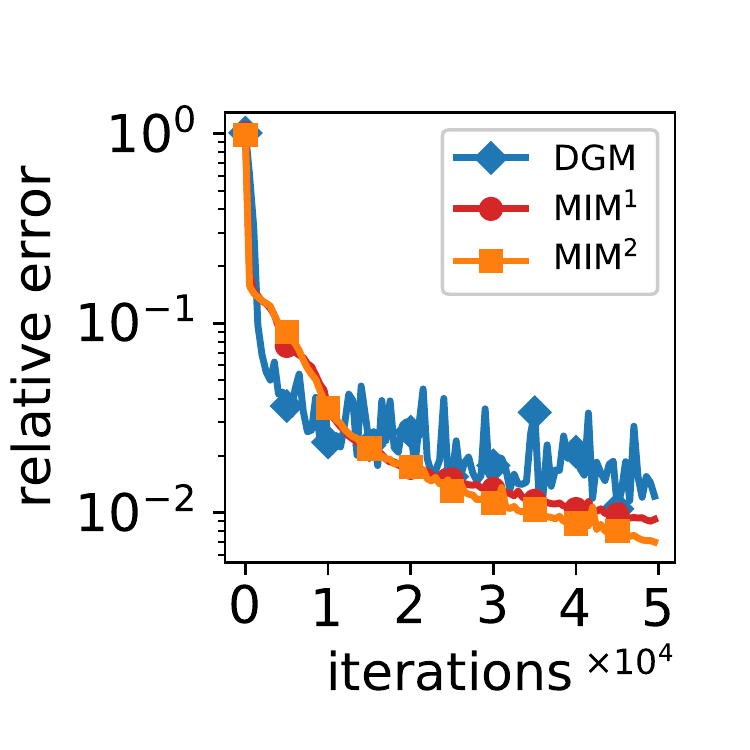}
		}
		\subfigure[16D]{
			\includegraphics[width=0.3\textwidth]{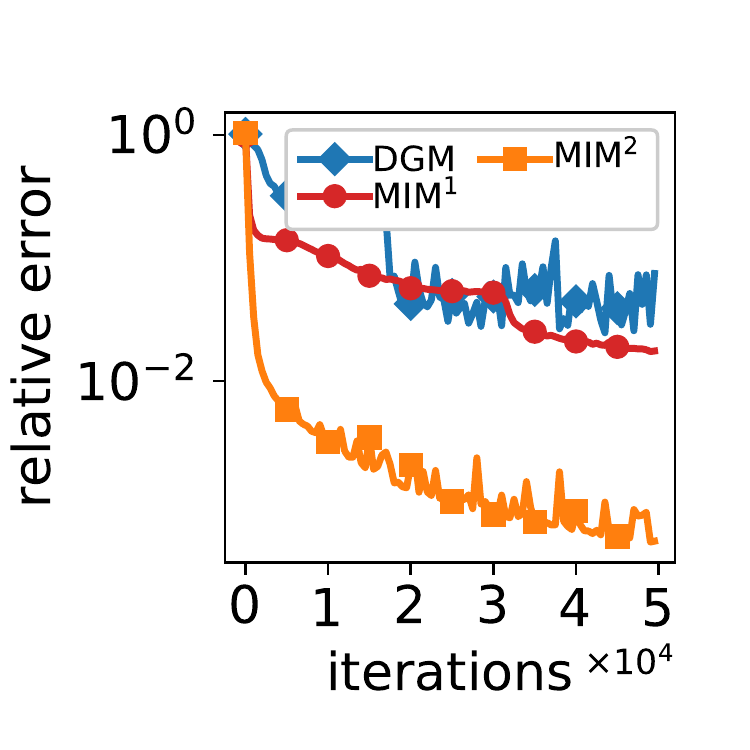}
		}
		\caption{Relative $L^2$ error of $u$ in terms of iteration number for Poisson equation defined in \eqref{eq:case_Neumann Boundary condition}.}
		\label{fig:NBCL2u}
	\end{figure}

\begin{figure}[ht]
	\centering
%
	\subfigure[4D]{
		\includegraphics[width=0.3\textwidth]{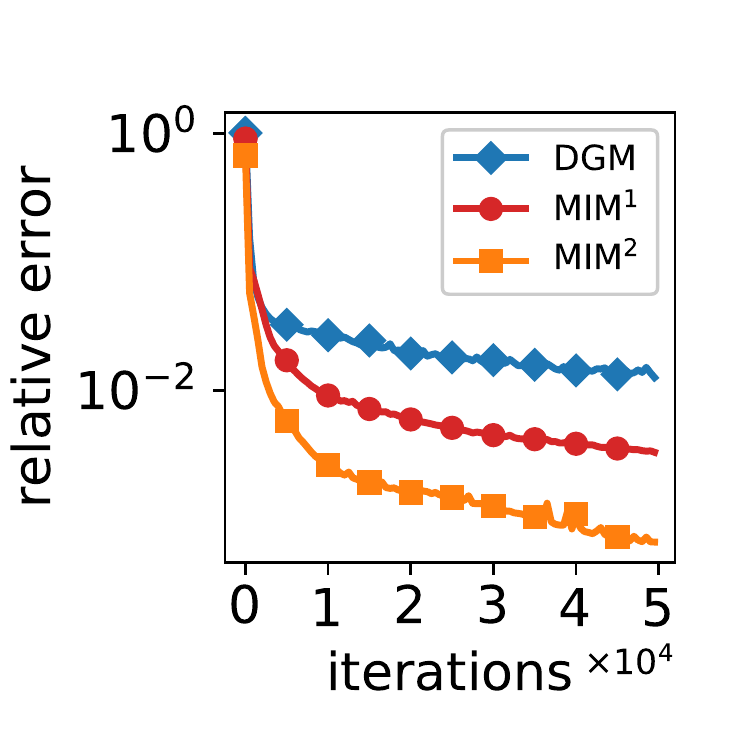}
	}
	\subfigure[8D]{
		\includegraphics[width=0.3\textwidth]{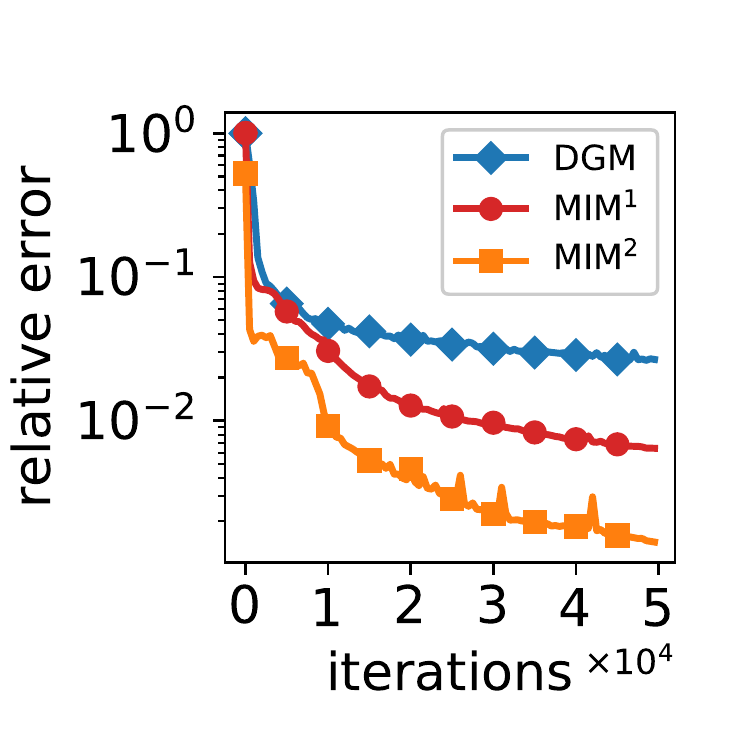}
	}
	\subfigure[16D]{
		\includegraphics[width=0.3\textwidth]{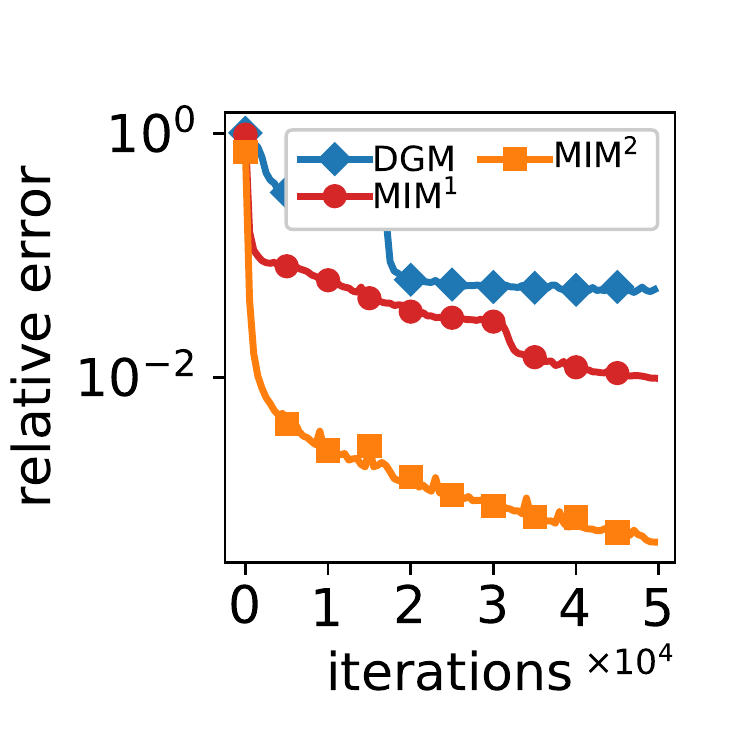}
	}
	\caption{Relative $L^2$ error of $\nabla u$ in terms of iteration number for Poisson equation defined in \eqref{eq:case_Neumann Boundary condition}.}
	\label{fig:NBCL2du}
\end{figure}

	For \eqref{eq:case_Neumann Boundary condition}, average errors of $u$ and $\nabla u$ over the last $100$ iterations are recorded in Table \ref{tbl:error of NBC}. The network depth $m=2$ and the activation function $x^2$ is used. Network widths are $5,10,15,20$ for $2,4,8,16$ dimensional problems, respectively. Time is recorded as the average CPU time per iteration. It is not surprising that MIM$^1$ costs less time than DGM since the DNN approximation in MIM satisfies the Neumann boundary condition automatically and both methods have similar network structures. It is surprising that MIM$^2$ costs less time than DGM since the number of parameters in MIM$^2$ is about twice of that in DGM. In terms of execution time, MIM$^1 < $ MIM$^2 < $ DGM.
	\begin{table}[ht]
		\centering
		\begin{tabular}{|c|c|c|c|c|}
			\hline
			\multirow{2}*{d} & \multirow{2}*{Method}& \multicolumn{2}{c|}{Relative $L^2$ error ($\times10^{-2}$)} &  \multirow{2}*{Time (s)} \\
			\cline{3-4}
			~& ~ &$u$ & $\nabla u$ & \\
			\hline
			\multirow{3}*{2}
			  & DGM &     0.3676 & 0.3714 & 0.04374\\
			~ & MIM$^1$ & 0.2941 & 0.1639 & 0.02925 \\
			~ & MIM$^2$ & 0.0565 & 0.0236 & 0.03514 \\
			\hline
			\multirow{3}*{4}
			  & DGM     & 1.0022  & 1.3272  & 0.07455\\
			~ & MIM$^1$ & 0.3751  & 0.3290  & 0.03603 \\
			~ & MIM$^2$ & 0.2294  & 0.0690  & 0.04141\\
			\hline
			\multirow{3}*{8}
			   & DGM & 2.0022 & 2.6551 & 0.13081\\
			~  & MIM$^1$ & 0.9049 & 0.6423 & 0.06642 \\
			~  & MIM$^2$ & 0.7261 & 0.1499 & 0.08716\\
			\hline
			\multirow{3}*{16} &DGM& 3.9796 & 5.0803 & 0.25621 \\
			~  & MIM$^1$ & 1.7631 & 1.0041 & 0.11082\\
			~  & MIM$^2$ & 0.0787 & 0.0236 & 0.15125\\
			\hline
		\end{tabular}
		\caption{Relative errors for $u$ and $\nabla u$ in DGM and MIM for Poisson equation defined in \eqref{eq:case_Neumann Boundary condition}.}\label{tbl:error of NBC}
	\end{table}
	Figure \ref{fig:NBCL2u} and Figure \ref{fig:NBCL2du} plot training processes of DGM and MIM in terms of relative $L^2$ errors for $u$ and $\nabla u $. Generally speaking, in terms of approximation error, MIM$^2 < $ MIM$^1 < $ DGM as expected. Therefore, MIM provides a better strategy over DGM. MIM provides better approximations in terms of relative $L^2$ errors for both $u$ and $\nabla u$. For $\nabla u$, the improvement of MIM$^1$ over DGM is about several times and that of MIM$^2$ over MIM$^1$ is about one order of magnitude. For $u$, the improvement is about several times. Moreover, a dimensional dependence is observed for both $u$ and $\nabla u$. The higher the dimension is, the better the approximation is.

	Table \ref{tbl:indep of D and A} records approximation errors of MIM and DGM in terms of activation function and network depth when $d=4$. MIM provides better approximations for both $\nabla u$ and $u$. It is not surprising that ReLU is not a suitable function for DGM due to high-order derivatives, but is suitable in MIM since only first-order derivatives are present in MIM.
	\begin{table}[ht]
		\centering\begin{tabular}{|c|c|c|c|c|c|c|c|}\hline
			\multirow{3}*{$\sigma$} & \multirow{3}*{$m$} &\multicolumn{6}{c|}{Relative $L^2$ error ($\times1$)} \\
			\cline{3-8}
			~&~  & \multicolumn{2}{c|}{DGM} & \multicolumn{2}{c|}{MIM$^1$} & \multicolumn{2}{c|}{MIM$^2$} \\
			\cline{3-8}
			~ & ~ & $u$ & $\nabla u$ & $u$ & $\nabla u$ & $u$ & $\nabla u$ \\
			\hline
			 \multirow{3}*{ReLU}
			    & 1 & 0.9197  & 0.9259 & 0.0890 & 0.0444  & 0.0264 & 0.0080 \\

			  ~ & 2 & 0.9210  & 0.9230 & 0.0245 & 0.0104  & 0.0265 & 0.0068\\

			  ~ & 3 & 0.9208  & 0.9216 & 0.0258 & 0.0113  & 0.0258 & 0.0084 \\
			  \hline
			 \multirow{3}*{ReQU}
			  & 1 & 0.0684 & 0.1003 & 0.0182 & 0.0127  & 0.0107 & 0.0042 \\

			~ & 2 & 0.0057 & 0.0118 & 0.0113 & 0.0047  & 0.0049 & 0.0017\\

			~ & 3 & 0.0124 & 0.0140 & 0.0040 & 0.0029  & 0.0042 & 0.0031\\
			\hline
			\multirow{3}*{ReCU}
			  & 1 & 0.4642 & 0.4644 & 0.0288 & 0.0159  & 0.0100 & 0.0033 \\

			~ & 2 & 0.0281 & 0.0170 & 0.0071 & 0.0055  & 0.0048 & 0.0013 \\

			~ & 3 & 0.0028 & 0.0031 & 0.0049 & 0.0036  & 0.0049 & 0.0013 \\
			\hline
		\end{tabular}
		\caption{Performance of MIM and DGM with respect to network depth and activation function for Poisson equation when $d=4$ . Network width is fixed to be $10$.}
		\label{tbl:indep of D and A}
	\end{table}

	\subsection{Monge-Amp\'ere equation}

	Consider the nonlinear Monge-Amp\'ere equation
	\begin{equation}\label{equ:MA}
		\left\{\begin{aligned}
		&\det(\nabla^2 u) = f(x) & x\in \Omega =  [-1,1]^d\\
		&u(x) = g(x) & x\in \partial \Omega
		\end{aligned} \right.
	\end{equation}
	with the exact solution defined as $u(x)= e^{1/d(\sum_{i=1}^d x_i^2)}$. Following Table \ref{tbl:loss DGM}, \ref{tbl:loss MIM} and \ref{tbl:loss BC}, we have the loss function in DGM
	\begin{equation*}
		L(u)= \|\det(\nabla^2u) - f\|_{2,\Omega}^2 + \lambda \|u-g\|^2_{2,\partial \Omega},
	\end{equation*}
	 and the loss function in MIM
	 \begin{equation*}
	 	L(u,p) = \|p - \nabla u \|_{2,\Omega}^2  + \|\det(\nabla p) - f\|_{2,\Omega}^2 + \lambda \|u-g\|^2_{2,\partial \Omega},
	 \end{equation*}
	 respectively. For \eqref{equ:MA}, the Dirichlet boundary condition can be enforced for both DGM and MIM. For comparison purpose, instead, we have the penalty term in both DGM and MIM. However, imposing exact boundary conditions is always encouraged in practice.
	
\begin{figure}[ht]
 	\centering
 	\subfigure[2D]{
 		\includegraphics[width=0.3\textwidth]{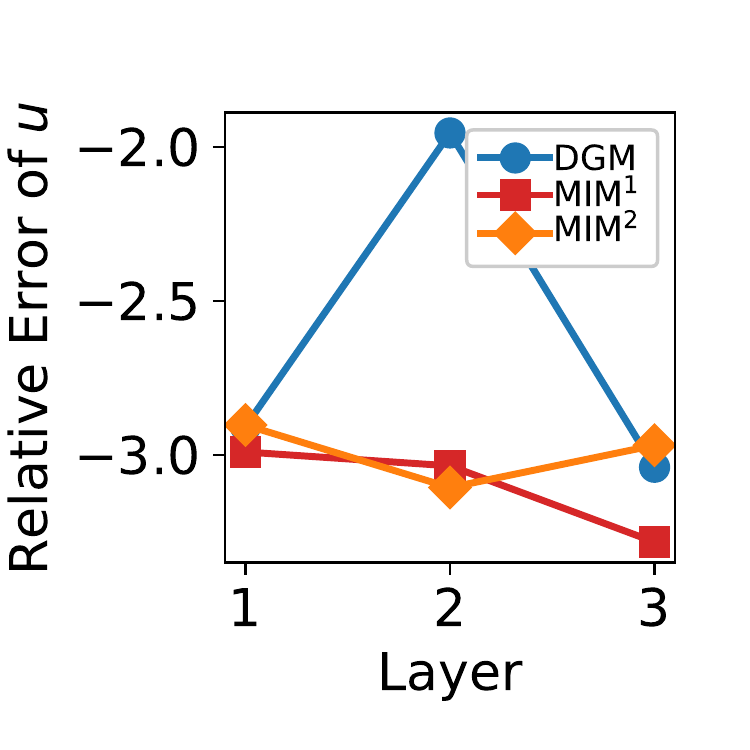}
 	}
 	\subfigure[4D]{
 		\includegraphics[width=0.3\textwidth]{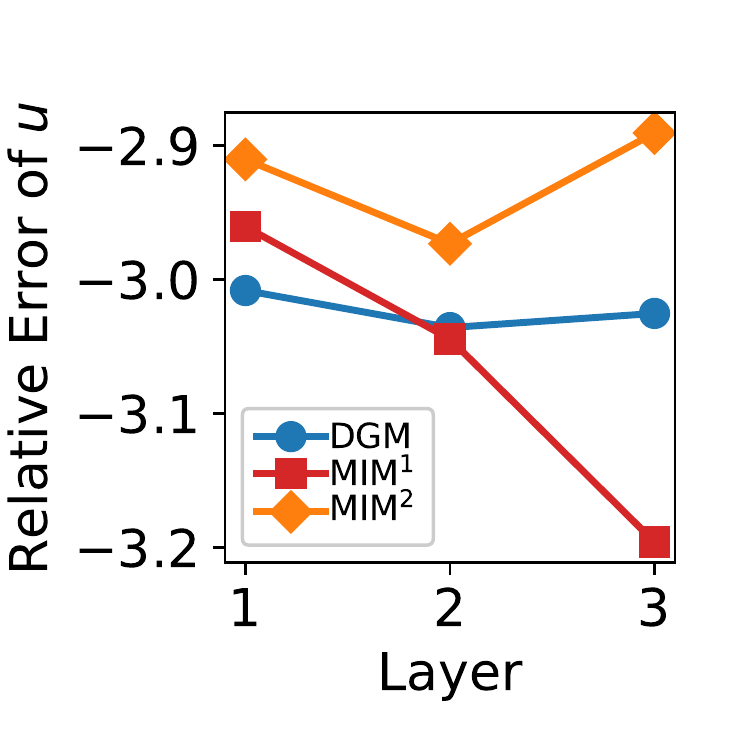}
 	}
 	\subfigure[8D]{
 		\includegraphics[width=0.3\textwidth]{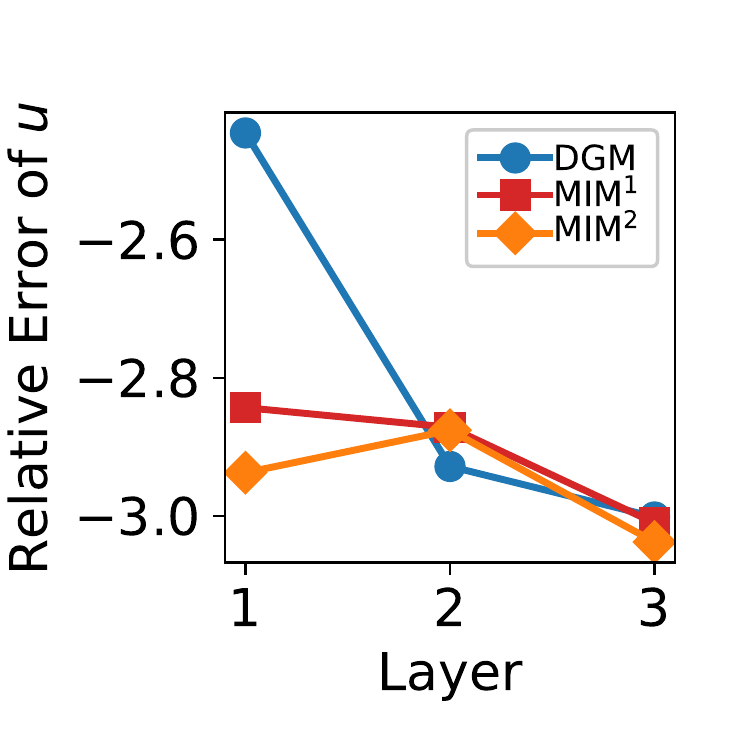}
 	}
 	\subfigure[2D]{
 		\includegraphics[width=0.3\textwidth]{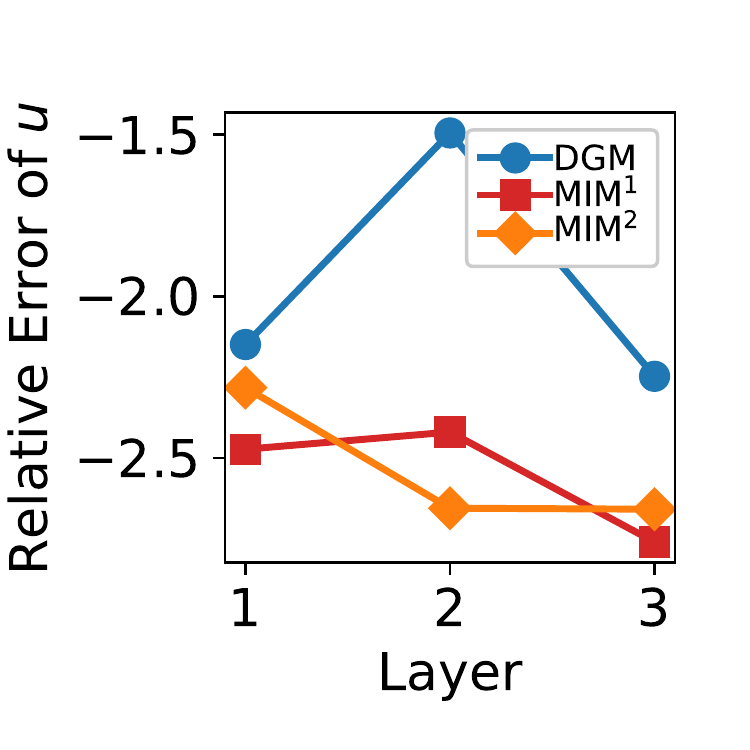}
 	}
 	\subfigure[4D]{
 		\includegraphics[width=0.3\textwidth]{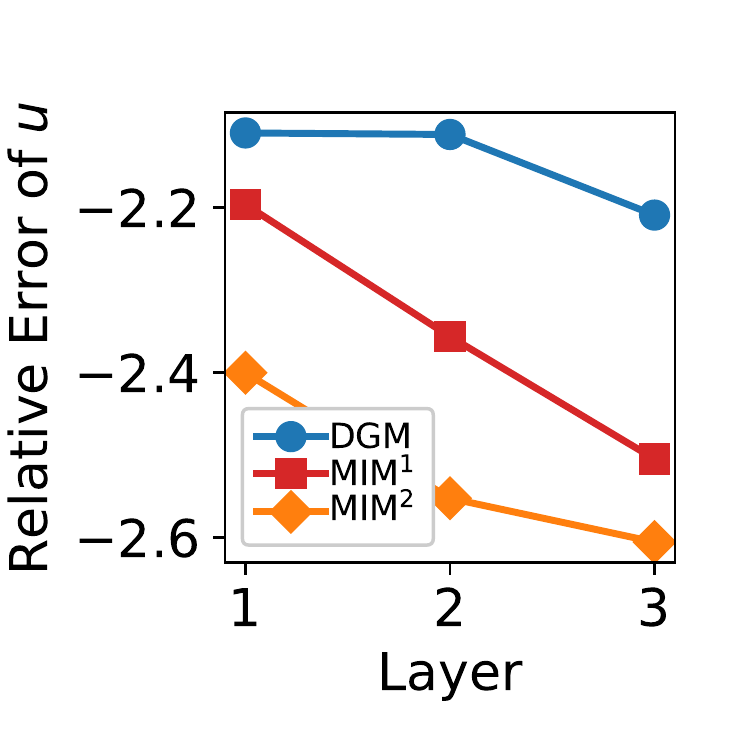}
 	}
 	\subfigure[8D]{
 		\includegraphics[width=0.3\textwidth]{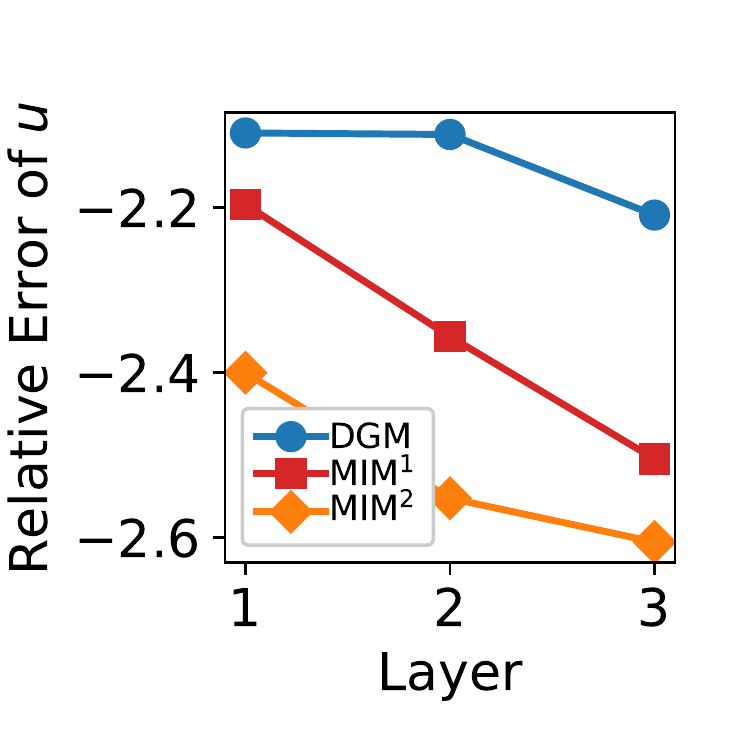}
 	}
 	\caption{Relative $L^2$ errors of $u$ and $\nabla u$ for Monge-Amp\'ere equation defined in \eqref{eq:case_Neumann Boundary condition}.}
 	\label{fig:indep of W and D in MA}
 \end{figure}

	 In this example, we fix the network depth $m = 2$ and the activation function as $\sigma(x) = \mathrm{ReQU}(x)$. Relative $L^2$ errors in the last $1000$ iterations with respect to the network width in different dimensions are recorded in Table \ref{tab:MA}. Figure \ref{fig:indep of W and D in MA} plots errors in terms of network width for different dimensions. The advantage of MIM is obvious from these results.
	 \begin{table}
	 	\centering\begin{tabular}{|c|c|c|c|c|c|c|c|}\hline
	 		\multirow{3}*{d} & \multirow{3}*{$n$} &\multicolumn{6}{c|}{Relative $L^2$ error ($\times 10^{-2}$)} \\
	 		\cline{3-8}
	 		~&~  & \multicolumn{2}{c|}{DGM} & \multicolumn{2}{c|}{MIM$^1$} & \multicolumn{2}{c|}{MIM$^2$} \\
	 		\cline{3-8}
	 		~ & ~ & $u$ & $\nabla u$ & $u$ & $\nabla u$& $u$ & $\nabla u$ \\
	 		\hline
	 		\multirow{3}*{2}
	 		 & 10 & 0.1236 & 0.7430 & 0.1023 & 0.3433& 0.1251 & 0.5218 \\	
	 		~& 20 & 1.1100 & 3.1940 & 0.0922 & 0.3804& 0.0784 & 0.0221\\
	 		~& 30 & 0.0913 & 0.5656 & 0.0522 & 0.1740& 0.1075 & 0.0219 \\
	 		\hline
	 		\multirow{3}*{4}
	 		    & 20 &  0.0981 & 0.7764 & 0.1095 & 0.6359& 0.1230 & 0.3977 \\
	 		~   & 30 &  0.0921 & 0.7731 & 0.0903 & 0.4399& 0.1063 & 0.2802\\
	 		~   & 40 &  0.0943 & 0.6174 & 0.0636 & 0.3127& 0.1287 & 0.2480\\
	 		\hline
	 		\multirow{3}*{8}
	 		& 30 & 0.3584 & 3.3902 & 0.1435 & 1.6318 & 0.1155 & 0.5170 \\
	 		~	& 40 & 0.1179 & 1.4663 & 0.1344 & 1.0721 & 0.1330 & 0.4873\\
	 		~	& 50 & 0.0997 & 1.2483 & 0.0977 & 0.8289 & 0.0917 & 0.4174\\
			\hline
	 	\end{tabular}
	 	\caption{Relative $L^2$ errors in the last $1000$ iterations with respect to the network width for Monge-Amp\'ere equation defined in \eqref{equ:MA}
	    for different dimensions. The network depth is fixed to be $m = 2$ and the activation function is fixed to be $\sigma(x) = \mathrm{ReQU}(x)$. }
	 	\label{tab:MA}
	 \end{table}

	\subsection{Biharmonic equation}
	
	Consider the biharmonic equation
	\begin{equation}\label{equ:biharmonic equation}
	\left\{\begin{aligned}
	&\Delta^2 u = \frac{\pi^4}{16} \sum_{k=1}^{d} \sin (\frac{\pi}{2}x) & x\in\Omega\\
	& u(x) = \sum_{k=1}^{d}\sin(\frac{\pi x}{2}) & x\in\partial \Omega \\
	& \frac{\partial u}{\partial n} = 0 & x\in\partial \Omega
	\end{aligned}\right.
	\end{equation}
	with the exact solution $u(x) = \sum_{k=1}^{d}\sin(\frac{\pi x}{2})$ over $\Omega = [-1,1]^d$.
	The loss function in DGM is
		\begin{equation*}
	L(u) = \|\Delta^2 u - \frac{\pi^4}{16} \sum_{k=1}^{d} \sin (\frac{\pi}{2}x) \|_{2,\Omega}^2 + \lambda_1 \| u - \sum_{k=1}^{d}\sin(\frac{\pi x}{2})\|_{2,\partial \Omega}^2 + \lambda_2 \|  \frac{\partial u}{\partial n} \|_{2,\partial \Omega}^2.
	\end{equation*}
	The loss function in MIM$_a$ is
	\begin{multline}
	\label{equ:loss for biharmonic scheme1}
		L(u,p,q,w) = \|p - \nabla u \|_{2,\Omega}^2 + \|q - \nabla \cdot p \|_{2,\Omega}^2   + \|w - \nabla q \|_{2,\Omega}^2 \\
			+ \| \nabla \cdot w - \frac{\pi^4}{16} \sum_{k=1}^{d} \sin (\frac{\pi}{2}x) \|_{2,\Omega}^2
			 + \lambda_1\| u - \sum_{k=1}^{d}\sin(\frac{\pi x}{2})\|_{2,\partial \Omega}^2 + \lambda_2\| p \|_{2,\partial \Omega}^2 ,
	\end{multline}
	and the loss function in MIM$_p$ is
	\begin{multline}
	\label{equ:loss for biharmonic scheme2}
	L(u,q) = \|q- \Delta u\|_{2,\Omega}^2 + \| \Delta  q - \frac{\pi^4}{16} \sum_{k=1}^{d} \sin (\frac{\pi}{2}x) \|_{2,\Omega}^2 \\
    + \lambda_1\| u - \sum_{k=1}^{d}\sin(\frac{\pi x}{2})\|_{2,\partial \Omega}^2 + \lambda_2\|  \frac{\partial u}{\partial n} \|_{2,\partial \Omega}^2.
	\end{multline}
	Again, we can enforce the exact boundary condition in MIM but cannot enforce it in DGM. For comparison purpose, we use penalty terms in both methods.
	
\begin{figure}
 	\centering
 	\subfigure[$u$]{
 		\includegraphics[width=0.45\textwidth]{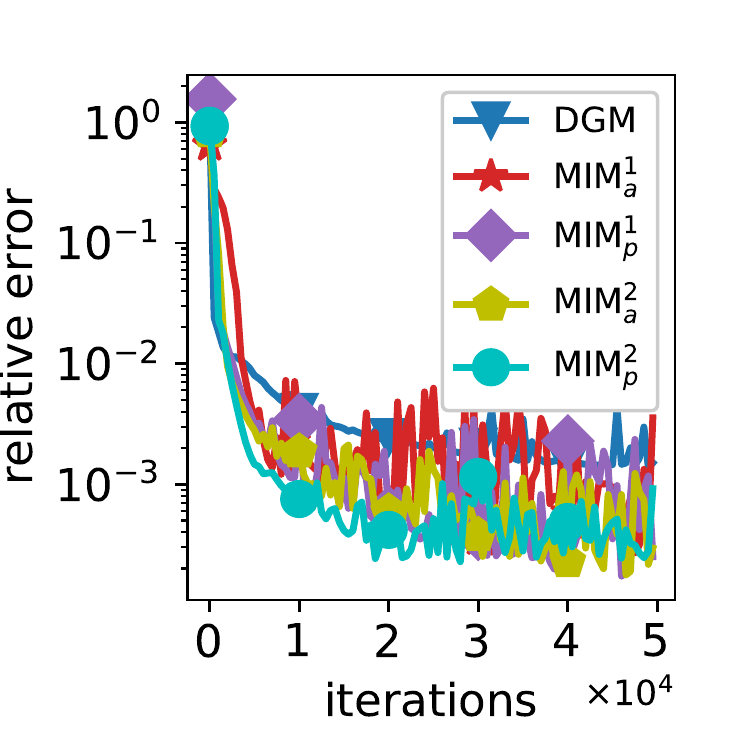}
 	}
 	\subfigure[$\nabla u$]{
	\includegraphics[width=0.45\textwidth]{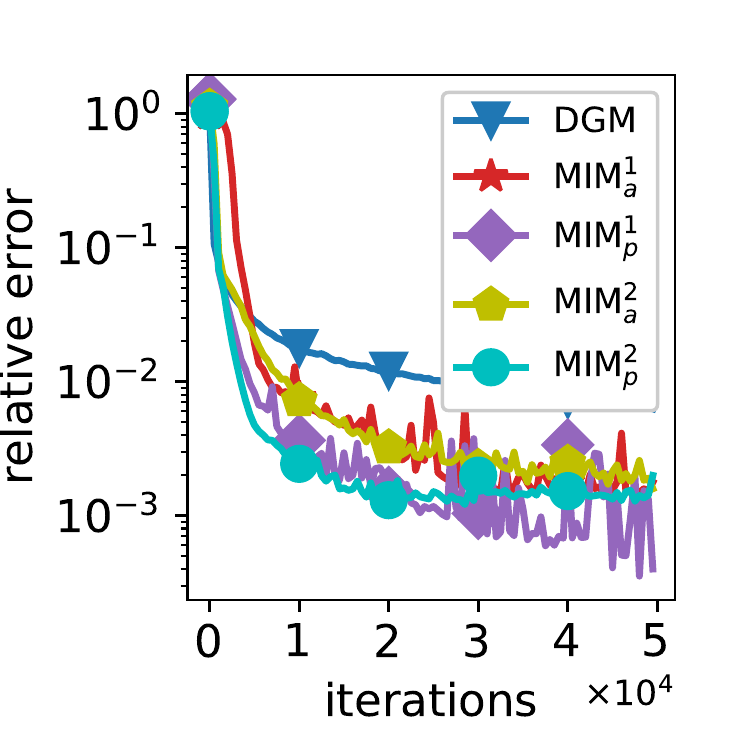}
	}
 	\subfigure[$\Delta u$]{
	\includegraphics[width=0.45\textwidth]{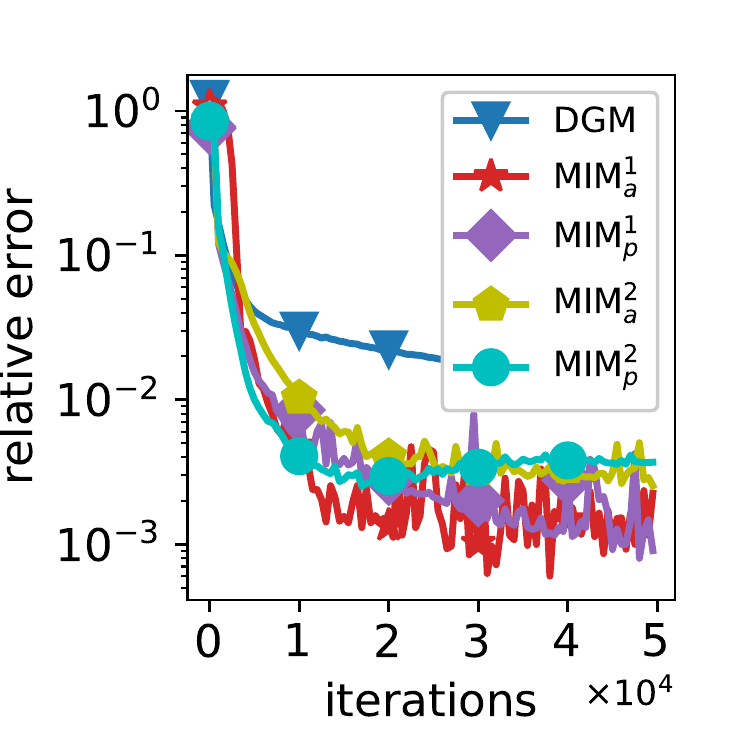}
	}
 	\subfigure[$\nabla \Delta u$]{
	\includegraphics[width=0.45\textwidth]{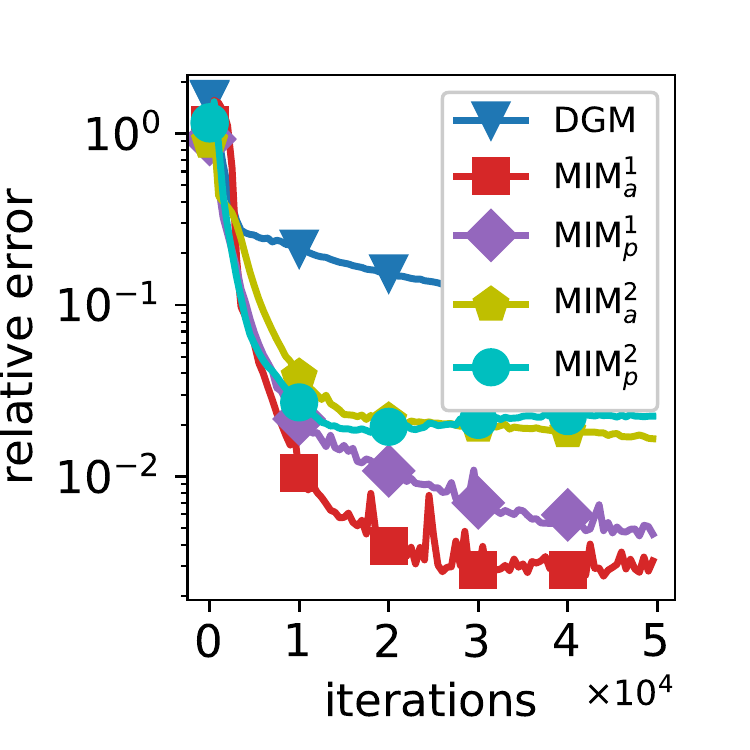}
	}
	\caption{Relative $L^2$ errors of $u$, $\nabla u$, $\Delta u$, $\nabla (\Delta u)$ in terms of iteration number for biharmonic equation. Both the solution and its derivatives are approximated by the same network in MIM$^1$, while different networks are used for the solution and its derivatives in MIM$^2$. MIM$_a$ means all derivatives are approximated and MIM$_p$ means only a subsect of derivatives ($\Delta u$ here) are approximated.}
	\label{fig:biharmonic 2D}
	\end{figure}

	Set $m=2$ and $n=8,10,20$ when $d=2,4,8$, respectively. Table \ref{tbl:error of biharmonic} records averaged errors in the last 1000 iterations.
	\begin{table}[ht]
		\centering
		\begin{tabular}{|c|c|c|c|c|c|c|}
			\hline
			\multirow{2}*{d} & \multirow{2}*{Method} & \multicolumn{4}{c|}{Relative $L^2$ error ($\times10^{-2}$ )} &  \multirow{2}*{Time (s)} \\
			\cline{3-6}
			~& ~ & $u$  & $\nabla u$  & $\Delta u$ & $\nabla (\Delta u)$ &\\
			\hline
			\multirow{5}*{2}
			  & DGM & 0.1656 & 0.6454 & 1.2333 & 8.8001 & 0.1034\\
			~ & MIM$^1_a$  & 0.1501 & 0.1929 & 0.1564 & 0.3067 & 0.1219 \\
			~ & MIM$^1_p$ & 0.0769 & 0.1155 & 0.1504 & 0.4984  & 0.1636\\
			~ & MIM$^2_a$ & 0.0526 & 0.2066 & 0.2937 & 1.6821 &  0.1393\\
			~ & MIM$^2_p$ & 0.0424 & 0.1417 & 0.3625 & 2.2231 & 0.2164 \\
			\hline
			\multirow{5}*{4}
			  & DGM       & 0.1330 & 0.6454 & 1.2333 & 8.8008 & 0.3292\\
			~ & MIM$^1_a$ & 0.4117 & 0.1929 & 0.1563 & 0.3066 & 0.2784\\
			~ & MIM$^1_p$ & 0.0845 & 0.1155 & 0.1504 & 0.4984 & 0.4692\\
			~ & MIM$^2_a$ & 0.1039 & 0.2066 & 0.2937 & 1.6821 & 0.2883\\
			~ & MIM$^2_p$ &	0.1111 & 0.1417	& 0.3625 & 2.2301 & 0.5919\\
			\hline
			\multirow{5}*{8}
			  & DGM     & 0.2488 & 1.0514 & 1.4594 & 13.4003 & 0.3292\\
			~ & MIM$^1_a$ & 0.3719 & 2.3855 & 0.6797 & 3.1015 & 0.2784\\
			~ & MIM$^1_p$ & 0.1856 & 0.6909 & 0.7840 & 4.7209 & 0.4692\\
			~ & MIM$^2_a$ & 0.1475 & 1.6657 & 1.2922 & 6.9594 & 0.8051\\
			~ & MIM$^2_p$ & 0.2881 & 0.9223 & 0.9981 & 6.4658 & 6.5148\\
			\hline
		\end{tabular}
		\caption{Relative errors for biharmonic equation defined in \eqref{equ:biharmonic equation}. MIM$_a$ and MIM$_b$ represent MIM with loss functions defined in \eqref{equ:loss for biharmonic scheme1} and \eqref{equ:loss for biharmonic scheme2}, respectively.}\label{tbl:error of biharmonic}
	\end{table}
	Relative $L^2$ errors for $u$ , $\nabla u$, $\Delta u$ and $\nabla (\Delta u)$ in terms of iteration number are plotted in Figure \ref{fig:biharmonic 2D} when $d=2$.
	Generally speaking, MIM provides better approximations for $u$, $\nabla u$, $\Delta u$, and $\nabla (\Delta u)$ than DGM. For MIM$_a$ and MIM$_p$, MIM$_p$ has a slightly better approximation accuracy comparable to that of MIM$_a$, although MIM$_a$ has $2d+2$ more outputs. These results are of interests since they are connected with results of local discontinuous Galerkin method that the formulation with a subset of derivatives has a better numerical performance \cite{Yan2002Dec,cockburn2009a}. We point out that MIM$_a$ has the advantage that the exact boundary condition can be enforced, although we use penalty terms for this example.
	
	\subsection{KdV equation}
	
	Consider a time-dependent linear KdV-type equation
	\begin{equation}\label{eq:kdv}
	\left\{
	\begin{aligned}
	&u_t + \sum_{k=1}^d u_{x_kx_kx_k} =0
	& (t,x)\in[0,T]\times\Omega\\
	&u(0,x) = u_0(x) = \sin(\sum_{k=1}^d x_k)
	& (t,x)\in[0]\times\Omega\\
	&u(t,x) \mathrm{\; is\; periodic\; in\; } x \\
	\end{aligned}\right.
	\end{equation}
	defined over $\Omega=[0,2\pi]^d$, where the exact solution $u(t,x)=\sin(\sum_{k=1}^d x_k+d t)$.
	We first rewrite it into the first-order system
	\begin{equation*}
	\begin{aligned}
	&p = \nabla u,\\
	&q = \mathrm{diag}(\nabla p ),\\	
	&u_t + \nabla \cdot q = 0.
	\end{aligned}
	\end{equation*}
    The loss function in DGM is
	\begin{equation*}
	\begin{aligned}
	L(u) &= \|u_t+\sum_{k=1}^du_{x_k x_k x_k}\|^2_{2,[0,1]\times\Omega} + \lambda_1 \|u-\sin(\sum_{k=1}^d x_k+d t)\|^2_{2,[0,1]\times\partial\Omega} \\\
	&  \quad + \lambda_2 \left(\sum_{k=1}^d\|u(x,t)-u(x\pm 2\pi e_k,t)\|^2_{2,\Omega}\right) \\
	&  \quad + \lambda_3 \left(\sum_{k=1}^d\|\nabla u(x,t)-\nabla u(x\pm2\pi e_k,t)\|^2_{2,\Omega}\right).
	\end{aligned}
	\end{equation*}
	Here $\{e_k\}_{k=1}^d$ is the standard basis set of $\mathbb{R}^d$. The loss function in MIM is
	\begin{equation*}
	\begin{aligned}
	L(u,p,q) & = \|p - \nabla u\|^2_{2,[0,1]\times\Omega} + \|q- \mathrm{diag}(\nabla p)\|^2_{2,[0,1]\times\Omega} \\
	& \quad + \|u_t+\nabla \cdot q\|^2_{2,[0,1]\times\Omega} + \lambda_1 \|u-\sin(\sum_{k=1}^d x_k+d t)\|^2_{2,[0,1]\times\partial\Omega} \\
	&  \quad + \lambda_2 \left(\sum_{k=1}^d\|u(x,t)-u(x\pm2\pi e_k,t)\|^2_{2,\Omega}\right) \\
	& \quad + \lambda_3\left(\sum_{k=1}^d\|p(x,t)-p(x\pm2\pi e_k,t)\|^2_{2,\Omega}\right).
	\end{aligned}
	\end{equation*}
	Relative $L^2$ errors of $u$, $\nabla u$, and $\mathrm{diag}(\nabla^2 u)$ are recorded in Table \ref{tbl:error of kdv}. Again, as shown in previous examples, MIM provides better results compared to DGM, especially for ReQU activation function. No obvious improvement of MIM$^2$ over MIM$^1$ is observed.
	\begin{table}[ht]
	\centering
	\begin{tabular}{|c|c|c|c|c|c|}
		\hline
		\multirow{2}*{$d$} & \multirow{2}*{$\sigma$} & \multirow{2}*{Method} & \multicolumn{3}{c|}{Relative $L^2$ error ($\times 10^{-2}$) }    \\
		\cline{4-6}
		~& ~ & ~ & $u$  & $\nabla u$  & $\mathrm{diag}(\nabla^2 u)$  \\
		\hline
		\multirow{6}*{1}
		& \multirow{3}*{ReQU}   
		& DGM     & 34.9171& 20.6788 & 34.3661 \\
		~ & ~ & MIM$^1$ & 0.5705 & 5.3709  & 0.5369 \\
		~ & ~ & MIM$^2$ & 1.2920 & 0.8129 &  1.9244 \\
		\cline{2-6}
		~ & \multirow{3}*{ReCU} 
		& DGM     & 0.7603 & 0.4785 & 0.5977 \\
		~ & ~ & MIM$^1$ & 0.0991 & 0.7313 & 0.0128 \\
		~ & ~ & MIM$^2$ & 0.5035 & 0.5804 & 0.1229 \\
		\hline
		\multirow{6}*{2}
		&\multirow{3}*{ReQU}  
		& DGM     & 84.8708 & 85.8114 & 85.8954 \\
		~ & ~ & MIM$^1$ & 2.9393  & 1.9996  & 2.9443 \\
		~ & ~ & MIM$^2$ & 2.1820  & 2.5591  & 2.1383 \\
		\cline{2-6}
		~ & \multirow{3}*{ReCU} &
		DGM     & 2.5483  & 2.1856 & 2.4431 \\
		~ & ~ & MIM$^1$ & 1.5410  & 2.3865 & 1.5645 \\
		~ & ~ & MIM$^2$ & 5.5900  & 5.7440 & 5.8957 \\
		\hline
		\multirow{6}*{3}
		& \multirow{3}*{ReQU} & DGM & 168.1755 & 168.1697 & 169.3528 \\
		& ~ & MIM$^1$ & 4.0421 & 4.0987 & 3.8496 \\
		& ~& MIM$^2$ & 7.7027 & 8.8787 & 9.1058 \\
		\cline{2-6}
		& \multirow{3}*{ReCU} & DGM & 1.9132 & 1.4846 & 1.7970 \\
		& ~ & MIM$^1$ & 1.5410 & 2.3865 & 1.5645 \\
		& ~ & MIM$^2$ & 5.5900 & 5.7440 & 5.8957 \\
		\hline
	\end{tabular}
	\caption{Relative $L^2$ errors for KdV equation defined in \eqref{eq:kdv}.}\label{tbl:error of kdv}
\end{table}

\section{Conclusion and Discussion}
\label{sec:conclusion}
	
	Motivated by classical numerical methods such as local discontinuous Galerkin method, mixed finite element method, and least-squares finite element method, we develop a deep mixed residual method to solve high-order PDEs in this paper. The deep mixed residual method inherits several advantages of classical numerical methods:
	\begin{itemize}
		\item Flexibility for the choice of loss function;
		\item Larger solution space with flexible choice of deep neural networks;
		\item Enforcement of exact boundary conditions;
		\item Better approximations of high-order derivations with almost the same cost.
	\end{itemize}
    Meanwhile, the deep mixed residual method also provides a better approximation for the PDE solution itself. These features make deep mixed residual method suitable for solving high-order PDEs in high dimensions.

    Boundary condition is another issue which is important for solving PDEs by DNNs. Enforcement of exact boundary conditions not only makes the training process easier, but also improves the approximation accuracy; see \cite{berg_unified_2018, chen2020bc} for examples. The deep mixed residual method has the potential for imposing exact boundary conditions such as Neumann boundary condition, mixed boundary condition, and Robin boundary condition. All these conditions cannot be enforced exactly in deep Galerkin method. This shall be investigated in a subsequent work \cite{lyu2020}.
    
    So far, in the deep mixed residual method, only experiences from classical numerical methods at the basic level are transferred into deep learning. We have seen its obvious advantages. To further improve the deep mixed residual method, we need to transfer our experiences from classical numerical analysis at a deeper level. For example, the choice of solution space relies heavily on the choice of residual in order to maximize the performance of least-squares finite element method \cite{Bochev2015}. Many other connections exist in discontinuous Galerkin method \cite{Cockburn2000} and mixed finite element method \cite{Boffi2013}. For examples, since only first-order derivatives appear in the deep mixed residual method, ReLU works well for all time-independent equations we have tested but does not work well for KdV equation. Therefore, it deserves a theoretical understanding of the proposed method in the language of linear finite element method \cite{he2018relu}. Another possible connection is to use the weak formulation of the mixed residual instead of least-squares loss, as done in deep learning by \cite{zang2020weak} and in discontinuous Galerkin method by \cite{Cockburn2000}. Realizing these connections in the deep mixed residual method will allow for a systematic way to understand and improve deep learning for solving PDEs.

\section{Acknowledgments}
	This work was supported by National Key R\&D Program of China (No. 2018YFB0204404) and National Natural Science Foundation of China via grant 11971021. We thank Qifeng Liao and Xiang Zhou for helpful discussions.

\bibliographystyle{amsplain}
\bibliography{ref}
\end{document}